\documentclass[review]{elsarticle}

\usepackage{lipsum}
\usepackage{amsfonts}
\usepackage{graphicx}
\usepackage{bm}
\usepackage{epstopdf}
\usepackage{fancyhdr}
\usepackage{xcolor}
\usepackage{tikz} 
\usetikzlibrary{shapes,patterns} 
\usepackage{algorithmic}
\usepackage[utf8]{inputenc}
\usepackage[linesnumbered,ruled,vlined]{algorithm2e} 
\SetKwInput{KwInput}{Input}    
\SetKwInput{KwOutput}{Output}  
\ifpdf%
\DeclareGraphicsExtensions{.eps,.pdf,.png,.jpg}
\else
\DeclareGraphicsExtensions{.eps}
\fi
\usepackage{amsopn}
\DeclareMathAlphabet{\mathpzc}{OT1}{pzc}{m}{it}
 
\usepackage{booktabs}
\usepackage{lineno,hyperref}
\modulolinenumbers[5]
\newdefinition{rmk}{Remark}
\DeclareMathOperator{\ddiv}{div}

\journal{Journal of \LaTeX\ Templates}

\begin{document}

\begin{frontmatter}

\title{A Local Fourier Analysis for Additive Schwarz Smoothers}

\author{\'Alvaro P\'e de la Riva \fnref{UZ}}
\fntext[UZ]{IUMA and Applied Mathematics Department, University of Zaragoza, Spain (apedelariva@gmail.com, carmenr@unizar.es, fjgaspar@unizar.es).}
\author{Carmen Rodrigo \fnref{UZ}} 
\author{Francisco J. Gaspar \fnref{UZ}} 
\author{James H. Adler \fnref{Tufts}} 
\fntext[Tufts]{Tufts University, Medford, Massachusetts, USA (james.adler@tufts.edu, xiaozhe.hu@tufts.edu).}
\author{Xiaozhe Hu \fnref{Tufts}} 
\author{Ludmil Zikatanov \fnref{PennState}} 
\fntext[PennState]{Penn State, University Park, Pennsylvania, USA (ludmil@psu.edu).}
%


\begin{abstract}
In this work, a local Fourier analysis is presented to study the convergence of multigrid methods based on additive Schwarz smoothers. This analysis is presented as a general framework which allows us to study these smoothers for any type of discretization and problem. 
The presented framework is crucial in practice since it allows one to know a priori the answer to questions such as what is the size of the patch to use within these relaxations, the size of the overlapping, or even the optimal values for the weights involved in the smoother.  Results are shown for a class of additive and restricted additive Schwarz relaxations used within a multigrid framework applied to high-order finite-element discretizations and saddle point problems, which are two of the contexts in which these type of relaxations are widely used. 
\end{abstract}


\end{frontmatter}


\section{Introduction}
\label{sec:intro}

Multigrid methods are among the best-known iterative solution techniques due to their demonstrated high efficiency for a wide range of problems. They accelerate the convergence of classical iterative methods by combining them with a coarse-grid correction technique. The design of efficient multigrid methods, however, depends crucially on the choice of their components. One of the most important ingredients of a multigrid algorithm is the so-called smoother or relaxation procedure, which often consists of a classical iterative method such as Jacobi or Gauss-Seidel.

Within a multigrid framework, a natural extension of point-wise smoothers are patch-wise smoothers. In order to apply such a relaxation, the computational domain is divided into small (overlapping or non-overlapping) patches, and then, one smoothing step consists of solving local  problems on each patch one-by-one either in a Jacobi-type or Gauss-Seidel-type manner. This results in an additive or multiplicative Schwarz smoother, respectively. 
One of the best-known multigrid methods based on this type of relaxation was proposed by Vanka in \cite{vanka} for solving the steady-state incompressible Navier-Stokes equations in primitive variables, discretized by a finite-volume scheme on a staggered grid. The computational domain is divided into cells with pressure nodes at the cell centers and velocity nodes at the cell faces. The smoothing procedure is a so-called symmetric coupled Gauss-Seidel technique (SCGS), which consists of solving local problems for each cell involving all the unknowns located at the cell. This is done cell by cell in a Gauss-Seidel-type manner and, therefore, can be viewed as a multiplicative Schwarz iteration.

Additive Schwarz-type iteration methods have been studied in \cite{Schoberl} as smoothers in a multigrid method for saddle point problems. It is shown that, under suitable conditions, the iteration can be interpreted as a symmetric inexact Uzawa method.
Restrictive additive Schwarz methods (RAS) were introduced  in \cite{Cai} as an efficient alternative to the classical additive Schwarz preconditioners. Convergence of RAS methods was proven in \cite{Frommer}, where it was shown that this method reduces communication time while maintaining the most desirable properties of the classical Schwarz methods. RAS preconditioners are widely used in practice and are implemented in several software packages. In \cite{Saberi}, a restrictive Schwarz method was proposed as smoother for solving the Stokes equations. It was observed that this smoother achieves comparable convergence rates to the multiplicative version, while being computationally less expensive per iteration.  In general, the class of additive and restricted additive Schwarz smoothers are characterized by their ability to deal with high-order discretizations, saddle point problems, and equations where the terms grad-div or curl-curl dominate.  Thus, in this work, we aim to study such relaxations within a multigrid framework.

Local Fourier analysis (LFA), or local mode analysis, is a commonly used approach for analyzing the convergence properties of geometric multigrid methods. In this analysis an infinite regular grid is considered and boundary conditions are not taken into account. LFA was introduced by Brandt in \cite{Bra77}, and afterwards extended in \cite{Bra94}. A good introduction can be found in the paper by St\"uben and Trottenberg \cite{Stuben} and in the books by Wesseling \cite{wesseling2004}, Trottenberg et al. \cite{trottenberg2001}, and Wienands and Joppich \cite{Wie01}.  It is the main
quantitative analysis for the convergence of multilevel algorithms, and results in a very useful tool for the design of multigrid methods. Moreover, in \cite{JCM-37-340} it has been recently proved that under standard assumptions LFA is a rigorous analysis, providing the exact asymptotic convergence factors of the multigrid method. 

LFA for multiplicative Schwarz smoothers was first performed in \cite{Siva} for the staggered finite-difference discretization of the Stokes equations, and in \cite{Molenaar1991} for a mixed finite-element discretization of the Laplace equation. In \cite{Rodrigo2016}, an LFA for multiplicative Schwarz smoothers on triangular grids is presented, and in  \cite{MacLachlan}, the analysis for such overlapping block smoothers is performed on rectangular grids for finite-element discretizations of the grad-div, curl-curl and Stokes equations. 

Whereas LFA has been widely applied to multiplicative Schwarz relaxations, few works on additive Schwarz smoothers can be found in the literature.  Whereas the Fourier modes are eigenfunctions of the multiplicative Schwarz smoothers (see \cite{MacLachlan} for a rigorous proof), this statement is not true for additive Schwarz smoothers.
In \cite{Boonen2008LocalFA}, an LFA to analyze a multicolored version of an additive Schwarz smoother for a curl-curl model problem was proposed. A non-standard LFA  to analyze this type of smoothers for the Stokes equations with P2-P1 and Q2-Q1 discretizations is considered in \cite{Farrel2021}. This analysis was also used in \cite{Greif2023} in order to study an additive Vanka-type smoother within a multigrid framework for the Poisson equation discretized by a standard finite-difference scheme. The analysis developed in both references can be seen as a particular case of the general framework analysis presented in this work, which allows us to study this class of smoothers for any type of discretization and problem.

Finally, we note that a non-standard LFA technique to predict the convergence rate of multigrid solvers for problems involving random and jumping coefficients was proposed in \cite{Prashant}.  The novelty of this new approach lies in the use of a specific basis of the Fourier space, rather than the standard basis which is based on the Fourier modes. This is the approach that we consider in this work and, as it will be shown, it allows us to propose an LFA to study the convergence of multigrid methods based on additive Schwarz smoothers for any discretization and problem. This general framework is crucial in practice, since it can answer questions such as what the appropriate size of the patch to use is, how to determine the size of the overlap, what the optimal values for the weights involved in the smoother are, and if the restricted additive version is the best choice. 

The rest of the paper is organized as follows. In Section \ref{sec:additive}, the additive and restricted additive Schwarz smoothers are introduced together with their corresponding algorithms. Section \ref{sec:lfa} is devoted to presenting the basis for the LFA performed in this work. Sections \ref{sec:high} and \ref{sec:saddle} deal with the application of the proposed local Fourier analaysis to high-order finite-element discretizations for scalar problems and saddle point type problems, respectively.  Finally, a summary and concluding remarks are given in Section \ref{sec:conclusions}.

\section{Additive Schwarz methods}
\label{sec:additive}
In this section we introduce the additive Schwarz methods used to solve a linear system of $n$ algebraic equations, $A u = b$. Consider a decomposition of the unknowns into subsets $\mathcal{B}_i$, $i = 1, \ldots,  nb$, such that each unknown, $u_j$, in vector $u$ is included in at least one block $\mathcal{B}_i$. For each subset of unknowns, $\mathcal{B}_i$, let $V_i: R^n \rightarrow R^{n_i}$ be the mapping from the global vector of unknowns to the ones given in $\mathcal{B}_i$, where $n_i$ is the size of $\mathcal{B}_i$. Here, $V_i$ is represented by a rectangular $n_i \times n$ matrix and, with a slight abuse of notation, we also refer to this matrix as $V_i.$  Its transpose, $V_i^T$, is a prolongation operator from $R^{n_i}$ to $R^n$ by padding it with zeros. Then,  $A_i=V_i A V_i^T$ is the local matrix corresponding to the unknowns in $\mathcal{B}_i$, that is, the restriction of $A$ to $\mathcal{B}_i$. For each block $\mathcal{B}_i$, the following small local system has to be solved, 
$$
A_i \delta {u}_i^j = V_i(b-A {u}^j),
$$
where $\delta {u}_i^j$ denotes the corrections' vector for the unknowns involved in block $\mathcal{B}_i$. 

These corrections can be computed via a multiplicative or an additive Schwarz method. Although multiplicative Schwarz methods provide better asymptotic convergence factors, their additive counterparts are more efficient since they are parallelizable. This is due to the fact that the corrections provided by solving the other local systems are not required to compute the corrections associated with a given block. The algorithm of the additive Schwarz method method is given as follows: \\

\begin{algorithm}[H]\label{alg:1}
\DontPrintSemicolon
  
  \KwInput{$u^j$. \textbf{Output:} $u^{j+1}$.} 
  $ r = b - A u^j$ \\
  \For{$i=1:nb$}
    {
        \textbf{Solve:} $A_i \delta u_i^j = V_i r.$ 
    }
    $  u^{j+1}=u^j + \sum_{i=1}^{nb} V_i^T D_i \delta u_i^j,$

\caption{Additive Schwarz method}
\end{algorithm}

\vspace{0.5cm}
\noindent In Algorithm \ref{alg:1}, $D_i$ is, in the general case, a diagonal weighting matrix. The diagonal entries $d_{ jj}$, $j = 1, \ldots, n$, of this matrix are usually taken as the inverse of the number of blocks sharing the unknown $u_j$, but can be tuned to improve the convergence factor of the method.

The restricted additive Schwarz method is based on the  use of prolongation operators, $\widetilde{V}_i^T$, in such a way that each entry $u_j$ of the vector unknown $u$ occurs in $\widetilde{V}_i u$  for exactly one $i$. This means that prolongation operators $\widetilde{V}_i^T$ are introduced in such a way that
$$
\sum_{i=1}^{nb} \widetilde{V}_i^T V_i = I,
$$
is satisfied. Then, the algorithm of the restricted additive Schwarz method is: \\

\begin{algorithm}[H]
\DontPrintSemicolon
  
  \KwInput{$u^j$. \textbf{Output:} $u^{j+1}$.} 
  $ r = b - A u^j$ \\
  \For{$i=1:nb$}
    {
        \textbf{Solve:} $A_i \delta u_i^j = V_i r.$ 
    }
    $  u^{j+1}=u^j + \sum_{i=1}^{nb} \widetilde{V}_i^T D_i \delta u_i^j,$

\caption{Restricted additive Schwarz method}
\end{algorithm}

\vspace{0.5cm}
When using these type of smoothers for solving a particular problem, some key questions arise. 
First, what should the size of the blocks be?  Second, how large should the overlap be?  Finally, the weights involved in the method are input parameters, and are extremely important to obtain an efficient and robust solver. Thus, determining the optimal values is a big concern.  In the next section, we present a general framework based on LFA that answers these questions, allowing us to know a priori the type of additive Schwarz method to use.

\section{Local Fourier Analysis}
\label{sec:lfa}

In this section, we describe LFA in a setting which allows us to estimate the multigrid convergence factors
by using additive Schwarz methods as smoothers. Standard LFA assumes that Fourier modes are eigenfunctions of the operators involved in the multigrid algorithm, such as the discrete operator, the inter-grid transfer operators, and the relaxation procedure. However, we mention three examples where this conventional analysis is not directly applicable. The first exception consists of problems with non-constant coefficients, as models involving random and jumping coefficients. Secondly, discretizations with different stencils at distinct grid-points, appearing for example in high-order finite-element methods also pose challenges to the classical approach. Finally, standard LFA is not applicable when Fourier modes are not eigenfunctions of the smoother considered in the multigrid algorithm.  Additive Schwarz methods fall in this latter category. 

In \cite{Prashant}, a non-standard LFA technique to predict the convergence rate of multigrid solvers for problems involving random and jumping coefficients was proposed. While the analysis in that work was applied only to problems with non-constant coefficients, it can be also applied to high-order finite-element discretizations. In fact, the special LFA performed in \cite{cuadraticos} to study multigrid methods for quadratic finite-element methods can be seen as a particular case of the analysis presented in \cite{Prashant}. In addition, as will be shown here, it can also be applied for analyzing multigrid methods based on additive Schwarz smoothers. Next, we briefly describe this analysis, and refer to \cite{Prashant} for more details. 

Given an infinite and regular grid, $\mathcal{G}_h$ with grid size $h$, Fourier modes are defined by $\varphi_h(x, \theta) = e^{\imath {\mathbf \theta} x}$, with $x \in \mathcal{G}_h$ and $\theta \in \Theta_h :=( \pi\/h,\pi/h]$. The Fourier modes span the so-called Fourier space:
\begin{equation}
    \mathcal{F}_h(\mathcal{G}_h) = \hbox{span} \lbrace \varphi_h(x,\theta)\, | \, \theta \in \Theta_h \rbrace.
\end{equation}

The main idea in the analysis presented in  \cite{Prashant} 
is to consider a specific basis of the Fourier space, rather than the standard basis which is based on the Fourier modes. For this purpose, we consider a splitting of $\mathcal{G}_h$ into $n^d$ infinite subgrids, where $d$ is the dimension of the problem. We then take a fixed window of size $\underbrace{n \times \ldots \times n}_d$ and its periodic extension. The size, $n$, of the window has to be chosen appropriately, which we give in the context of additive Schwarz smoothers. For the sake of simplicity, we consider the same size $n$ in each direction.  We note that the analysis still holds even if they are different. Once $n$ is fixed, for every $\bm{k} = (k_1, \ldots, k_d), \;k_1, \ldots, k_d=0,\ldots,n-1$, subgrid, $G_h^{\bm{k}}$, is defined as follows:
\begin{equation}\label{inf_subgrids}
    G_h^{\bm{k}} = \lbrace \bm{k}h + (l_1,\ldots,l_d)nh \hspace{0.1cm} | \hspace{0.1cm} l_1,\ldots,l_d \in \mathbb{Z}\rbrace.
\end{equation}
For each frequency $\bm\theta^{\bm{0}} \in \Theta_{nh} := (-\pi/nh,\pi/nh]^d$, we introduce the following grid functions:
\begin{equation}\label{grid_func}
    \psi_h^{\bm k} (\bm\theta^{\bm{0}},\bm x) = \varphi_h (\bm\theta^{\bm{0}},\bm x) \, \chi_{{G}_h^{\bm k}}. 
\end{equation}
As in \cite{Prashant}, the Fourier space spanned by these functions,
\begin{equation}\label{fourier_space}
    \mathcal{F}_h^{n^d} (\bm\theta^{\bm{0}}) = \operatorname{span}\lbrace \psi_h^{\bm k}(\bm\theta^{\bm{0}},\cdot),\; \bm{k} = (k_1, \ldots, k_d),\; k_1, \ldots, k_d = 0,\ldots,n-1\rbrace,
\end{equation}
is the same Fourier space generated by the standard Fourier modes,
\begin{equation}\label{fourier_space_2}
\operatorname{span}\lbrace \varphi_h(\bm\theta^{\bm{0}}_{\bm k},\cdot),\; \bm\theta^{\bm{0}}_{\bm k} = \bm\theta^{\bm{0}} + {\bm k} \frac{2\pi}{nh},\; \bm{k} = (k_1, \ldots, k_d),\; k_1, \ldots, k_d = 0,\ldots,n-1\rbrace.
\end{equation} 
Due to the relation between the grid-functions of the given new Fourier basis and the standard
Fourier modes, it can be shown, in many situations, that the two-grid operator,
\begin{equation}\label{two-grid-operator}
M_h = (I_h-P_{2h}^hA_{2h}^{-1}R_h^{2h}A_h)S_h^{\nu},
\end{equation} 
satisfies the following invariance property for any frequency $\bm\theta^{\bm{0}} \in \Theta_{nh}$,
$$
M_h :  \mathcal{F}_h^{n^d} (\bm\theta^{\bm{0}}) \rightarrow  \mathcal{F}_h^{n^d} (\bm\theta^{\bm{0}}).
$$
In (\ref{two-grid-operator}),  $P_{2h}^h$ and $R_h^{2h}$ are the prolongation and restriction operators, $A_h$ and $A_{2h}$ are the fine- and coarse-grid operators, $I_h$ is the identity, $S_h$ is the smoothing operator, and  $\nu$ is the number of smoothing steps. 

When the usual analysis can be applied, and a standard coarsening $H=2h$ is considered, the Fourier space has dimension $2^d$, yielding the so called $2h-$harmonics spaces. In this case, $n=2$ in (\ref{fourier_space}), the size of the window is $2^{d}$, and both analyses coincide. The advantage of using the grid-functions in (\ref{grid_func}) is that it allows one to study problems for which the analysis based on the Fourier modes is not applicable. For example, using $n=2$ one could study the multigrid convergence for quadratic finite-element methods. As we will see in the next section, one can study finite-element methods of arbitrary order just taking a bigger value of $n$. Additionally, we will show that Additive Schwarz methods can be analyzed by considering  the grid-functions in (\ref{grid_func}), and choosing an adequate window size, i.e., an adequate value of $n$.

\section{High-order finite-element methods for scalar problems}
\label{sec:high}
As a model problem, we consider the Poisson equation in a $d-$dimensional domain, $\Omega = (0,1)^d$, with homogeneous Dirichlet boundary conditions,
\begin{equation}\label{model_problem}
\begin{array}{ccc}
-\Delta u  &=& f, \ \  \mbox{in} \quad \Omega, \\ 
u &=& 0, \ \  \mbox{on} \quad \partial \Omega.
\end{array}
\end{equation}
The variational formulation of this problem reads as follows: Find $u \in H_0^1(\Omega)$ such that 
$$
a(u,v) = (f,v), \quad \forall v \in H_0^1(\Omega),
$$
where
$$
a(u,v) = \int_{\Omega} \nabla u \cdot \nabla v \, {\rm d} x, \quad\mbox{and}\quad 
(f,v) = \int_\Omega f v \, {\rm d} x.
$$
Let ${\cal T}_h$ be a partition of the domain $\Omega \subset R^d$, and associate to ${\cal T}_h$ the high-order finite-element space, $V_{hp} \subset H_0^1(\Omega)$ defined as $V_{hp} = \{ u_h \in H_0^1(\Omega) \;  | \;  u_h |_T \in {\mathbb Q}_p,\, \forall \, T \, \in {\cal T}_h\}$, where ${\mathbb Q}_p$ is the space of polynomials up to total degree $p$ on each variable. Thus, the Galerkin approximation of the variational problem is given by: Find $u_h \in V_{hp}$, such that 
\begin{equation} \label{approxvariational}
a(u_h,v_h) = (f,v_h), \quad \forall v_h \in V_{hp}.
\end{equation}
The solution of the Galerkin approximation problem is a linear combination, $u_h = \sum_{i=1}^{n_h} u_i \varphi_i$, where $\{\varphi_1, \ldots, \varphi_{n_h}\}$ is a basis for $V_{hp}$ and ${\rm dim}\, V_{hp} = n_h$. In order to compute ${\mathbf u} = (u_1,\ldots, u_{n_h})^T$, the following linear system must be solved:
$
A {\mathbf u} = {\mathbf b},
$
where $A = (a_{i,j}) = (a(\varphi_j,\varphi_i))_{i,j=1}^{n_h}$ is the so-called stiffness matrix and ${\mathbf b}  = (b_i) = (f,\varphi_i)_{i=1}^{n_h}$ is the right-hand side vector.

The goal now is to demonstrate the utility of the analysis presented in the previous section for the design of efficient geometric multigrid methods based on additive Schwarz relaxation schemes. For this purpose, we consider an $h$-multigrid, i.e., we keep  the polynomial order unchanged in the geometric mesh hierarchy, combined with the canonical high-order restriction and prolongation operators. There are, however, other possibilities to build a multilevel hierarchy for high-order discretized problems. Another option is to use a $p$-multigrid technique, which consists of constructing the coarse problems within the multigrid algorithm by reducing the polynomial degree of the finite-element space. This approach can be combined with a standard $h$-multigrid on the coarsest level $p=1$. Although we only present results here for the $h$-multigrid case, the analysis proposed in this work can be applied to any other alternative for high-order discretizations, as those approaches mentioned above.


\subsection{One-dimensional case}
First, we consider the simplest case when $d = 1$, the computational domain is the interval, $\Omega  = (0, 1)$, and the corresponding two-point boundary value problem is $- u''(x)  =f(x), \  x \in \Omega, \; u(0) = u(1) = 0$. We consider linear finite-element methods, i.e., $p=1$, and perform LFA based on infinite subgrids to estimate the convergence factor of the multigrid method using additive Schwarz methods as the relaxation. In additive Schwarz, we vary the size of the block, $k$, from $2$ to $7$ and the size of the overlap, $ov$, from $ov =1$ to the maximum overlapping $ov=k-1$. The diagonal elements of the weighting matrix, $D_i$, in Algorithm \ref{alg:1} are taken to be the natural weights, i.e., the reciprocal of the number of blocks in which each unknown appears. For example, if the size of the block is $k=3$ and the size of the overlap is $ov = 1$, the weights for the three unknowns of a block are $1/2$, $1$, $1/2$. In Table \ref{p1_1D_pesofijo}, we show the two-grid asymptotic convergence factors provided by LFA with only one smoothing iteration, $\nu=1$, for different block-sizes and overlap sizes, considering a linear finite-element discretization of the 1D Poisson equation. We do not show the experimentally computed asymptotic convergence factors because they are exactly the same as those provided by LFA. 
\begin{table}[htb!]
	\begin{center}		
		\begin{tabular}{ccccccc}
			\cline{2-7}
			   & \multicolumn{6}{c}{Block-size (k)}\\  
			\cline{1-7}
			 Overlap (ov)  & 2 & 3 & 4 & 5 & 6 & 7 \\
			\hline 
			$1$ & 0.33 & 0.99 & 0.40 & 0.50 & 0.43 & 0.50 \\
			$2$ & -      & 0.33 & 0.20 & 0.50 & 0.29 & 0.25 \\
			$3$ & -      & -      & 0.20 & 0.50 & 0.21 & 0.99 \\
			$4$ & -      & -      & -      & 0.20 & 0.14 & 0.25 \\
			$5$ & -      & -      & -      & -      & 0.14& 0.33 \\
			$6$ & -      & -      & -      & -      & -      & 0.14 \\ 
			\hline
		\end{tabular}	
\caption{Linear finite-element discretization, $p=1$, for the 1D Poisson problem. Two-grid asymptotic convergence factors applying $\nu=1$ smoothing steps of an additive Schwarz method.  Block-size, $(k)$, ranges from $2$ to $7$ and the size of the overlap, $ov$, from $1$ to $k-1$. Natural weights are taken in the smoother.}		
	\label{p1_1D_pesofijo}				
	\end{center}
\end{table}

From the results presented in Table \ref{p1_1D_pesofijo}, we observe that a bigger block-size or overlapping among blocks does not always provide a better asymptotic convergence factor.

Next, we consider the restrictive additive Schwarz method as the smoother within the multigrid framework. Now, the natural weights will be one or zero. For example, if the size of the block is $k=3$ and the size of the overlap is $ov = 1$, the weights for the three unknowns of the block are $1$, $1$, $0$. In Table \ref{1Dp1_factors_restric}, we show the asymptotic convergence factors provided by our two-grid analysis by using one smoothing step of the restrictive additive Schwarz method with natural weights. Block-sizes ranging from $2$ to $7$ and all the possible overlapping among the blocks are considered. Again, the experimentally computed asymptotic convergence factors are the same as those provided by the LFA, and thus not included in the table. 
\begin{table}[htb!]
	\begin{center}
		\begin{tabular}{ccccccc}
			\cline{2-7}
			   & \multicolumn{6}{c}{Block-size}\\  
			\cline{1-7}
			 Overlapping & 2 & 3 & 4 & 5 & 6 & 7 \\
			\hline 
			$1$ & 0.75 & 1.00 & 0.40 & 0.50 & 0.43& 0.50 \\
			$2$ & -      & 1.00 & 0.60 & 0.50 & 0.57 & 0.44 \\
			$3$ & -      & -      & 0.87 & 0.87 & 0.28 & 1.00 \\
			$4$ & -      & -      & -      & 1.00 & 0.71 & 0.44 \\
			$5$ & -      & -      & -      & -      & 0.92 & 1.00 \\
			$6$ & -      & -      & -      & -      & -      & 1.00 \\ 
			\hline
		\end{tabular}
\caption{Linear finite-element discretization, $p=1$, for the 1D Poisson problem. Two-grid asymptotic convergence factors applying $\nu=1$ smoothing steps of a restrictive additive Schwarz method with natural weights. Block-size, $(k)$, ranges from $2$ to $7$ and the size of the overlap, ($ov$), from $1$ to $k-1$.}
			\label{1Dp1_factors_restric}
	\end{center}
\end{table}

The obtained results are not good in general, observing even no convergence of the multigrid mehod for some combinations of block-size and size of the overlap. The advantage of having the LFA tool provided in this work, however, is that this analysis can find the optimal weights to improve the performance of the multigrid method. Hence, in Table \ref{1Dp1_factors_restric_opt}, we show the asymptotic convergence factors together with the optimal weights provided by the analysis. Again, the two-grid analysis is based on only one smoothing step of the restrictive additive Schwarz methods with optimal weights, and different block-sizes and sizes of the overlap among the blocks are considered.
\begin{table}[htb!]
	\begin{center}
		\begin{tabular}{ccccccc}
			\cline{2-7}
			   & \multicolumn{6}{c}{Block-size}\\  
			\cline{1-7}
			 Ov & 2 & 3 & 4 & 5 & 6 & 7 \\
			\hline 
			$1$ & 0.45 (0.6) & 0.34(0.67) & 0.17 (0.83) & 0.20 (0.8) & 0.18 (0.82) & 0.20 (0.8) \\
			$2$ & -      &  0.37 (0.68) & 0.15 (0.71) & 0.20 (0.8) & 0.22 (0.78) & 0.18 (0.82) \\
			$3$ & -      & -      & 0.43 (0.71) & 0.34 (0.66) & 0.16 (0.84) & 0.34 (0.67) \\
			$4$ & -      & -      & -      &  0.36 (0.66) & 0.20 (0.7) & 0.18 (0.82) \\
			$5$ & -      & -      & -      & -      & 0.40 (0.7) & 0.34 (0.66) \\
			$6$ & -      & -      & -      & -      & -      & 0.34 (0.66) \\ 
			\hline
		\end{tabular}
\caption{Linear finite-element discretization, $p=1$, for the 1D Poisson problem. Two-grid asymptotic convergence factors (and optimal weights in parenthesis) applying one iteration of the restrictive additive Schwarz method with optimal weights. Block-size, $(k)$, ranges from $2$ to $7$ and the size of the overlap, ($ov$), from $1$ to $k-1$.}
			\label{1Dp1_factors_restric_opt}
	\end{center}
\end{table}

Comparing the results in Tables \ref{1Dp1_factors_restric} and \ref{1Dp1_factors_restric_opt},  we observe a significant improvement of the asymptotic convergence factors when appropriate non-trivial weights are used. Once again, better asymptotic convergence factors are not directly related to bigger block-sizes or larger overlapping among blocks.

Our analysis, though, allows us to extend LFA to high-order discretizations. Here, we fix the size of the block and the overlap and we vary the polynomial order in the finite-element space. An element-based additive Schwarz method is considered, where each block contains all the basis functions with support in the corresponding element. Therefore, the size of the blocks is chosen as $p+1$, where $p$ denotes the polynomial degree. The size of the overlap is fixed to be the minimum, i.e., $ov=1$. In Table \ref{1D_additive_element}, we show the asymptotic convergence factors estimated by the two-grid analysis by using one smoothing step of the additive Schwarz method (AS) and one step of the restricted additive Schwarz method (RAS). The polynomial degrees range from $p=2$ to $p=8$.
\begin{table}[htb!]
	\begin{center}
		\begin{tabular}{cccccccc}
			\cline{2-8} 
			& \multicolumn{7}{c}{Polynomial degree}\\  
			\cline{1-8}
			 Smoother & $p=2$ & $p=3$ & $p=4$ & $p=5$ & $p=6$ & $p=7$ & $p=8$\\
			\hline 
			AS & 0.50 & 0.09 & 0.17 & 0.11 & 0.17 & 0.12 & 0.16 \\ 
		        RAS & 0.50 & 0.12 & 0.15 & 0.11 & 0.17 & 0.12 & 0.16 \\ 
			\hline
		\end{tabular}
\caption{High-order finite-element discretizations for the 1D Poisson problem. Asymptotic convergence factors predicted by the two-grid analysis by using one smoothing step of the additive Schwarz method (AS) and the restricted additive Schwarz method (RAS) with natural weights. The size of the block is $k=p+1$. Minimum overlapping, $ov=1$, is considered and the polynomial degree ranges from $p=2$  to $p=8$.}
\label{1D_additive_element}
	\end{center}
\end{table}

For each polynomial degree, $p$, we find that the restricted additive Schwarz smoother achieves comparable convergence rates to the additive Schwarz smoother, having better properties in terms of scalability and applicability to high-performance computing. We conclude that the restricted additive Schwarz method is a favorable smoother for high-order discretizations of elliptic equations. The authors are not aware of any work comparing these approaches.

\begin{rmk}
The size of the window, $n$, in the Fourier analysis is the least common multiple of $2p$ and the minimum number $(k-ov)j$ such that $(k-ov)j > k$, with $j$ being any natural number. In the particular case that the size of the block $k=p+1$ and $ov =1$, the size of the block is the least common multiple of $2p$ and the minimum number $pj$ such that $pj > p+1$. For example, if $p=2$ the size of the windows is $n=4$.
\end{rmk}

\subsection{Two-dimensional case}
We now consider the two-dimensional case of model problem (\ref{model_problem}).  We follow the same element-based smoothing strategy using additive Schwarz methods. Thus, the size of local linear systems will be $k=(p+1)^2$ and the overlapping among the blocks in both directions is still one, $ov=1$. In Table \ref{2D_additive_element}, we show the asymptotic convergence factors predicted by the two-grid analysis by using $\nu=2$ smoothing steps of the corresponding  AS and the RAS methods with natural weights.  The polynomial degree ranges from $p=1$ to $p=8$.
\begin{table}[htb]
	\begin{center}
		\begin{tabular}{ccccccccc}
			\cline{2-9} 
			& \multicolumn{8}{c}{Polynomial degree}\\  
			\cline{1-9}
			 Smoother & $p=1$ & $p=2$ & $p=3$ & $p=4$ & $p=5$ & $p=6$ & $p=7$ & $p=8$\\
			\hline 
			AS & 0.14 & 0.15 & 0.21 & 0.25 & 0.31 & 0.36 & 0.40 & 0.43 \\ 
		        RAS  & 0.19 & 0.15 & 0.21 & 0.25 & 0.31 & 0.36 & 0.40 & 0.43 \\ 
			\hline
		\end{tabular}
\caption{High-order finite-element discretizations for the 2D Poisson problem. Asymptotic convergence factors predicted by the two-grid analysis using two smoothing steps of the additive Schwarz method (AS) and the restricted additive Schwarz method (RAS) with natural weights. The size of the blocks are $k=(p+1)^2$. Minimum overlapping $ov=1$ is assumed and the polynomial degree ranges from $p=1$  to $p=8$.}
\label{2D_additive_element}
	\end{center}
\end{table}

Comparing the results of both smoothers in Table \ref{2D_additive_element}, we conclude that the restrictive case is better from a computational cost point of view, as it avoids communication between processors. Next, we consider the performance of the multigrid method in terms of the number of smoothing steps at each level. In Table \ref{ciclos2D}, we show the two-grid asymptotic convergence factors provided by LFA,  $\rho_{2g}$, together with the asymptotic convergence factors computationally obtained by our own implementation, $\rho_{h}$, using $V(1,0)$, $V(1,1)$, $V(2,1)$ and $V(2,2)$ cycles. The coarsest grid used here within the $V$-cycle consists of only one interior grid point. Again, we employ the element-based restrictive additive Schwarz smoother with natural weights for polynomial degree ranging from $p=1$ to $p=8$.  We see that the asymptotic convergence factors predicted by the LFA match with high accuracy those computationally obtained by applying $V$-cycles.  

\begin{table}[h]
	\begin{center}
		\begin{tabular}{ccccccccc}
			\cline{2-9}  
			 & \multicolumn{2}{c}{$V(1,0)$} & \multicolumn{2}{c}{$V(1,1)$} & \multicolumn{2}{c}{$V(2,1)$} & \multicolumn{2}{c}{$V(2,2)$}\\
			 \hline 
			  $p$ & $\rho_{2g}$ & $\rho_{h}$ & $\rho_{2g}$ & $\rho_{h}$ & $\rho_{2g}$ & $\rho_{h}$ & $\rho_{2g}$ & $\rho_{h}$ \\ 
			\hline 
			$1$ & 0.41 & 0.41 & 0.19 & 0.18 & 0.10 & 0.09 & 0.03 & 0.03 \\ 
			$2$ & 0.39 & 0.40 & 0.15 & 0.15 & 0.06 & 0.06 & 0.02 & 0.02 \\ 
			$3$ & 0.46 & 0.45 & 0.21 & 0.21 & 0.09 & 0.09 & 0.04 & 0.04 \\ 
			$4$ & 0.50 & 0.50 & 0.25 & 0.25 & 0.13 & 0.13 & 0.06 & 0.06 \\ 
			$5$ & 0.56 & 0.56 & 0.31 & 0.31 & 0.18 & 0.18 & 0.10 & 0.10 \\ 
			$6$ & 0.60 & 0.60 & 0.36 & 0.36 & 0.21 & 0.21 & 0.13 & 0.13 \\ 
			$7$ & 0.63 & 0.63 & 0.40 & 0.40 & 0.26 & 0.25 & 0.16 & 0.16 \\ 
			$8$ & 0.66 & 0.66 & 0.43 & 0.43 & 0.28 & 0.28 & 0.19 & 0.19 \\ 
			\hline
		\end{tabular}	
\caption{High-order finite-element discretizations for the 2D Poisson problem. Two-grid asymptotic convergence factors provided by LFA, $\rho_{2g}$, together with the computed asymptotic convergence factors, $\rho_{h}$, using $V(1,0)$, $V(1,1)$, $V(2,1)$, $V(2,2)$ cycles with the element-based restrictive additive Schwarz methods and natural weights. The polynomial degree ranges from $p=1$ to $p=8$.}	
			\label{ciclos2D}			
	\end{center}
\end{table}

\section{Saddle point problems}
\label{sec:saddle}
In practice, Schwarz methods are among the most commonly used smoothers for solving saddle point problems. Whereas LFA has been widely applied to multiplicative Schwarz relaxations, few works on additive Schwarz smoothers are found in the literature. In \cite{Farrel2021}, LFA was used to analyze this type of smoothers for the Stokes equations discretized by Taylor-Hood P2-P1 and Q2-Q1 schemes. This analysis can be seen as a particular case of the one introduced here, when using a window of size $2\times2$. 
In this section, we show that the presented Fourier analysis can be also applied to study the convergence of multigrid methods based on additive Schwarz methods applied to such saddle-point problems.  

For this purpose, we consider the quasi-static Biot's model for poroelasticity.  This model assumes  that we have a deformable porous medium, which is linearly elastic, isotropic, homogeneous, and saturated by an incompressible Newtonian fluid. Given these assumptions, the well-known displacement-pressure formulation \cite{wang} must satisfy the following equations in a domain, $\Omega \subset \mathbb{R}^d$, 
\begin{equation}
\begin{array}{ccc}
& & -\ddiv\ \boldsymbol \sigma'  +   \alpha \, { \nabla \,  p} = \rho {\mathbf g}, \qquad
\boldsymbol \sigma' = 2\mu \boldsymbol \varepsilon(\mathbf{u})  + \lambda\ddiv(\mathbf{u}), \\
& & \displaystyle \frac{\partial}{\partial t} \left(\frac{1}{M}  p + \alpha { \nabla \cdot \, {\textbf{u}}} \right)- \nabla \cdot \left( \frac{1}{\mu_f} {\bm K}  (\nabla p - \rho_f {\mathbf g}) \right)= f, 
\end{array}
\end{equation} 
where $\boldsymbol \sigma'$ and $\boldsymbol \varepsilon$ are the effective stress and strain
tensors, $\lambda$ and $\mu$ are the Lam\'e coefficients, $p$ is the pore pressure, $\mathbf{u}$ is the displacement vector field, ${\mathbf g}$ is the gravity tensor, $\alpha$ is the Biot-Willis constant (which we assume is equal to one), ${\bm K}$ is the permeability of the porous medium, $\mu_f$ is the fluid viscosity, $M$ is the Biot modulus, and $f$ is a source term.
To complete the formulation of the problem, we add appropriate boundary and initial conditions. For instance,
\begin{equation}\label{bound-cond}
\begin{array}{ccccc}
  p = 0, & &  \quad \boldsymbol \sigma' \, {\bm n} &=& {\bm t}, \quad \hbox {on }\Gamma _t, \\
{\bm u} = {\bm 0}, & & \quad  \displaystyle  {\bm K}  \left(\nabla p  - \rho_f {\bm g} \right) \cdot {\bm n} &=& 0, \quad \hbox {on } \Gamma _c,
\end{array}
\end{equation}
where ${\bm n}$ is the unit outward normal to the boundary and
$\Gamma_t \cup \Gamma_c = \Gamma:=\partial\Omega$, with $\Gamma_t$ and $\Gamma_c$
disjoint subsets of $\Gamma$ having non null measure.
For the initial time, $t=0$, the following condition is fulfilled,
\begin{equation}\label{ini-cond}
     \left (\frac{1}{M} p + \alpha \nabla \cdot {\bm u} \right)\, (\bm{x},0)=0, \,  \bm{x} \in\Omega.
\end{equation}

 To discretize the problem, we use Taylor-Hood finite elements in space, ${\mathbb Q}_2- {\mathbb Q}_1$, and a backward-Euler scheme in time. Let ${\cal T}_h$ be a triangulation of $\Omega$ composed of rectangles. Defining the discrete spaces as
\begin{eqnarray*}
{\bm V}_h &=& \{ {\bm u}_h \in (H^1(\Omega))^d \ |  \; \forall T \in {\cal T}_h,  \; {\bm u}_{h}|_T \in {\mathbb Q}_2^d, \ {\bm u}_h|_{\Gamma _c} = {\bm 0} \}, \\
Q_h &=& \{p_h \in H^1(\Omega) \ |  \; \forall T \in {\cal T}_h, \; p_{h}|_T \in {\mathbb Q}_1 \  p_h|_{\Gamma _t} = 0\},
\end{eqnarray*}
the fully discretized scheme at time $t_m, \; m = 1,2,\ldots$, is written as follows:
find $({\bm u_h^m} ,p_h^m) \in {\bm V}_h \times Q_h$ such that
\begin{eqnarray}
    & \hskip -90pt a(\bm{u}_h^m,\bm{v}_h) -  \alpha(p_h^m,\ddiv \bm{v}_h) = (\rho \bm g,\bm{v}_h),
     ~ \forall \  \bm{v}_h \in \bm V_h, \label{discrete_variational1}\\
   &\alpha(\ddiv {\bar{\partial}}_t{\bm{u}_h^m},q_h) + b(p_h^m,q_h)   = (f_h^m,q_h) +  ( {\bm K}\mu_f^{-1} \rho_f {\bm g},\nabla q_h ),
   ~ \forall \ q_h \in Q_h,\label{discrete_variational2}
\end{eqnarray}
where ${\bar{\partial}}_t{\bm{u}_h^m} := ({\bm{u}_h^m} - {\bm{u}_h^{m-1}})/\tau$, with $\tau$ the time discretization parameter, $(\cdot,\cdot)$ is the standard inner product in the space $L_2(\Omega)$, and the bilinear forms $a(\cdot,\cdot)$ and $b(\cdot,\cdot)$ are given as
\begin{eqnarray*}\label{bilinear}
a(\bm{u},\bm{v}) &=& 2 \mu \int_{\Omega}\boldsymbol \varepsilon(\bm{u}):\boldsymbol \varepsilon(\bm{v}) \, {\rm d} \Omega +
\lambda\int_{\Omega} \ddiv\bm{u}\ddiv\bm{v} \, {\rm d} \Omega, \\
b(p,q) &=&  \int_{\Omega} \frac{{\bm K}}{\mu_f} \nabla p \cdot \nabla q \, {\rm d} \Omega.
\end{eqnarray*} 
This fully discrete scheme leads to a saddle-point system of equations at each time step of the form 
${\cal A}x = b$, where matrix
${\cal A}$ is a $2\times2$ block symmetric indefinite matrix:
\begin{equation}\label{system}
{\cal A}=  \left(
\begin{array}{cc}
A & B^T \\
B & -C
\end{array}
\right),
\end{equation}
with matrices A and C both being symmetric and positive definite.

To solve this system of equations, we consider a geometric multigrid method with standard coarsening based on an additive Schwarz method as the relaxation scheme. For transfer operators we consider the canonical restriction and prolongation operators. The classical Schwarz smoother for this type of saddle-point problems is typically defined as a set of blocks that consists of one pressure unknown and all the velocity unknowns that are connected to it.  That is, the degrees of freedom corresponding to the nonzero entries in the ith row of $B$ plus the $i$-th pressure degree of freedom. This yields a local system of size $51 \times 51$ and $76 \times 76$ for two-dimensional and three-dimensional problems, respectively. Due to the high computational cost associated with this type of smoother, the paralellization of the relaxation method becomes crucial for real applications. Therefore, the additive and the restricted additive Schwarz smoothers are a more natural choice than multiplicative versions, and it is clear that a Fourier analysis tool for them is of great importance for the design of geometric multigrid methods for saddle point problems. 

\begin{figure}[htb]
	\begin{center} 
	 \begin{tikzpicture}[domain=0:4,scale=0.5]
      \draw[very thin,color=gray!50,step=0.8cm,] (0,0) grid (8.8,8.8); 
      \draw [ step=0.8cm,color=gray!70, thick] (1.5999,1.5999) grid (7.2,7.2);
      
      \draw [pattern=crosshatch,pattern color=black!80, thick](1.6,2.4) -- (7.2,2.4);
      \draw [pattern=crosshatch,pattern color=black!80, thick](1.6,4) -- (7.2,4);
      \draw [pattern=crosshatch,pattern color=black!80, thick](1.6,5.6) -- (7.2,5.6); 
      \draw [pattern=crosshatch,pattern color=black!80, thick](3.2,1.6) -- (3.2,7.2);
      \draw [pattern=crosshatch,pattern color=black!80, thick](4.8,1.6) -- (4.8,7.2);
      \draw [pattern=crosshatch,pattern color=black!80, thick](6.4,1.6) -- (6.4,7.2);
        
      \draw[fill=black!80,color=black] (1.6,7.2) circle (.1);
      \node[  align=center,thick] at   (2.4,7.2)  {\scriptsize{$\square$}}; 
      \draw[fill=black!80,color=black] (3.2,7.2) circle (.1);
      \node[  align=center,thick] at   (4,7.2)  {\scriptsize{$\square$}}; 
      \draw[fill=black!80,color=black] (4.8,7.2) circle (.1);
      \node[  align=center,thick] at   (5.6,7.2)  {\scriptsize{$\square$}}; 
      \draw[fill=black!80,color=black] (6.4,7.2) circle (.1);
      \node[  align=center,thick] at   (7.2,7.2)  {\scriptsize{$\square$}};  
      
      \node[  align=center,thick] at  (1.6,6.4) {\tiny{$\bigcirc$}}; 
      \node[  align=center,thick] at  (2.4,6.4) {$\times$}; 
      \node[  align=center,thick] at  (3.2,6.4) {\tiny{$\bigcirc$}}; 
      \node[  align=center,thick] at  (4,6.4) {$\times$}; 
      \node[  align=center,thick] at  (4.8,6.4) {\tiny{$\bigcirc$}}; 
      \node[  align=center,thick] at  (5.6,6.4) {$\times$}; 
      \node[  align=center,thick] at (6.4,6.4) {\tiny{$\bigcirc$}}; 
      \node[  align=center,thick] at  (7.2,6.4) {$\times$};

      \draw[fill=black!80,color=black] (1.6,5.6) circle (.1);
      \node[  align=center,thick] at   (2.4,5.6)  {\scriptsize{$\square$}}; 
      \draw[fill=black!80,color=black] (3.2,5.6) circle (.1);
      \node[  align=center,thick] at   (4,5.6)  {\scriptsize{$\square$}}; 
      \draw[fill=black!80,color=black] (4.8,5.6) circle (.1);
      \node[  align=center,thick] at   (5.6,5.6)  {\scriptsize{$\square$}}; 
      \draw[fill=black!80,color=black] (6.4,5.6) circle (.1);
      \node[  align=center,thick] at   (7.2,5.6)  {\scriptsize{$\square$}};  
      
      \node[  align=center,thick] at  (1.6,4.8) {\tiny{$\bigcirc$}}; 
      \node[  align=center,thick] at  (2.4,4.8) {$\times$}; 
      \node[  align=center,thick] at  (3.2,4.8) {\tiny{$\bigcirc$}}; 
      \node[  align=center,thick] at  (4,4.8) {$\times$}; 
      \node[  align=center,thick] at  (4.8,4.8) {\tiny{$\bigcirc$}}; 
      \node[  align=center,thick] at  (5.6,4.8) {$\times$}; 
      \node[  align=center,thick] at (6.4,4.8) {\tiny{$\bigcirc$}}; 
      \node[  align=center,thick] at  (7.2,4.8) {$\times$};

      \draw[fill=black!80,color=black] (1.6,4) circle (.1);
      \node[  align=center,thick] at   (2.4,4)  {\scriptsize{$\square$}}; 
      \draw[fill=black!80,color=black] (3.2,4) circle (.1);
      \node[  align=center,thick] at   (4,4)  {\scriptsize{$\square$}}; 
      \draw[fill=black!80,color=black] (4.8,4) circle (.1);
      \node[  align=center,thick] at   (5.6,4)  {\scriptsize{$\square$}}; 
      \draw[fill=black!80,color=black] (6.4,4) circle (.1);
      \node[  align=center,thick] at   (7.2,4)  {\scriptsize{$\square$}};  
      
      \node[  align=center,thick] at  (1.6,3.2) {\tiny{$\bigcirc$}}; 
      \node[  align=center,thick] at  (2.4,3.2) {$\times$}; 
      \node[  align=center,thick] at  (3.2,3.2) {\tiny{$\bigcirc$}}; 
      \node[  align=center,thick] at  (4,3.2) {$\times$}; 
      \node[  align=center,thick] at  (4.8,3.2) {\tiny{$\bigcirc$}}; 
      \node[  align=center,thick] at  (5.6,3.2) {$\times$}; 
      \node[  align=center,thick] at (6.4,3.2) {\tiny{$\bigcirc$}}; 
      \node[  align=center,thick] at  (7.2,3.2) {$\times$};

      \draw[fill=black!80,color=black] (1.6,2.4) circle (.1);
      \node[  align=center,thick] at   (2.4,2.4)  {\scriptsize{$\square$}}; 
      \draw[fill=black!80,color=black] (3.2,2.4) circle (.1);
      \node[  align=center,thick] at   (4,2.4)  {\scriptsize{$\square$}}; 
      \draw[fill=black!80,color=black] (4.8,2.4) circle (.1);
      \node[  align=center,thick] at   (5.6,2.4)  {\scriptsize{$\square$}}; 
      \draw[fill=black!80,color=black] (6.4,2.4) circle (.1);
      \node[  align=center,thick] at   (7.2,2.4)  {\scriptsize{$\square$}};  
      
      \node[  align=center,thick] at  (1.6,1.6) {\tiny{$\bigcirc$}}; 
      \node[  align=center,thick] at  (2.4,1.6) {$\times$}; 
      \node[  align=center,thick] at  (3.2,1.6) {\tiny{$\bigcirc$}}; 
      \node[  align=center,thick] at  (4,1.6) {$\times$}; 
      \node[  align=center,thick] at  (4.8,1.6) {\tiny{$\bigcirc$}}; 
      \node[  align=center,thick] at  (5.6,1.6) {$\times$}; 
      \node[  align=center,thick] at (6.4,1.6) {\tiny{$\bigcirc$}}; 
      \node[  align=center,thick] at  (7.2,1.6) {$\times$};

      \draw[color=black,very thick] (1.5999,1.5999) rectangle (7.2,7.2);

      \draw[fill=black!80,color=blue] (0,8.8) circle (.1);
        \node[  align=center,thick,blue] at (0.8,8.8) {\scriptsize{$\square$}};
        \draw[fill=black!80,color=blue] (1.6,8.8) circle (.1);
        \node[  align=center,thick,blue] at (2.4,8.8) {\scriptsize{$\square$}}; 
        \draw[fill=black!80,color=blue] (3.2,8.8) circle (.1);
        
        \node[  align=center,thick,blue] at (0,8) {\tiny{$\bigcirc$}};
        \node[  align=center,thick,blue] at (0.8,8) {$\times$};
        \node[  align=center,thick,blue] at (1.6,8) {\tiny{$\bigcirc$}};
        \node[  align=center,thick,blue] at (2.4,8) {$\times$}; 
        \node[  align=center,thick,blue] at (3.2,8) {\tiny{$\bigcirc$}};
        
        \draw[fill=black!80,color=blue] (0,7.2) circle (.1);
        \node[  align=center,thick,blue] at (0.8,7.2) {\scriptsize{$\square$}};
        \draw[fill=black!80,color=blue] (1.6,7.2) circle (.1);
        \node[  align=center,thick,blue] at (2.4,7.2) {\scriptsize{$\square$}}; 
        \draw[fill=black!80,color=blue] (3.2,7.2) circle (.1);
        
        \node[  align=center,thick,blue] at (0,6.4) {\tiny{$\bigcirc$}};
        \node[  align=center,thick,blue] at (0.8,6.4) {$\times$};
        \node[  align=center,thick,blue] at (1.6,6.4) {\tiny{$\bigcirc$}};
        \node[  align=center,thick,blue] at (2.4,6.4) {$\times$}; 
        \node[  align=center,thick,blue] at (3.2,6.4) {\tiny{$\bigcirc$}};
        
        \draw[fill=black!80,color=blue] (0,5.6) circle (.1);
        \node[  align=center,thick,blue] at (0.8,5.6) {\scriptsize{$\square$}};
        \draw[fill=black!80,color=blue] (1.6,5.6) circle (.1);
        \node[  align=center,thick,blue] at (2.4,5.6) {\scriptsize{$\square$}}; 
        \draw[fill=black!80,color=blue] (3.2,5.6) circle (.1); 
        
         \draw [ step=0.8cm,color=blue, thick] (-0.0001,5.5999) grid (3.2,8.8);
         \node[below] at (0,-0.25) {$7$}; 
         \node[below] at (0.8,-0.25) {$8$};
         \node[below ] at (1.6,-0.25) {$1$}; 
         \node[below ] at (2.4,-0.25) {$2$};
         \node[below ] at (3.2,-0.25) {$3$};
         \node[below ] at (4,-0.25) {$4$};
         \node[below ] at (4.8,-0.25) {$5$};
         \node[below ] at (5.6,-0.25) {$6$};
         \node[below ] at (6.4,-0.25) {$7$};
         \node[below ] at (7.2,-0.25) {$8$};
         \node[below ] at (8,-0.25) {$1$};
         
         \node[left ] at (-0.25,0.8) {$1$};
         \node[left ] at (-0.25,1.6) {$8$}; 
         \node[left ] at (-0.25,2.4) {$7$}; 
         \node[left ] at (-0.25,3.2) {$6$};
         \node[left ] at (-0.25,4) {$5$};
         \node[left ] at (-0.25,4.8) {$4$};
         \node[left ] at (-0.25,5.6) {$3$};
         \node[left ] at (-0.25,6.4) {$2$}; 
         \node[left ] at (-0.25,7.2) {$1$};
          \node[left ] at (-0.25,8) {$8$};
         \node[left ] at (-0.25,8.8) {$7$}; 
      \end{tikzpicture}
        \caption{Minimum size window for the two-dimensional Biot's model. The block centered on the upper-left pressure unknown of the grid is shown in blue.}
		\label{Stokes_window}
	\end{center}
\end{figure}
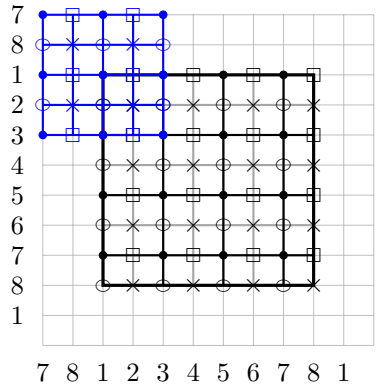

For two-dimensional problems, since the number of different stencils is $m=2$, and each basis function for the pressure shares common support with $5$ displacement unknowns on each direction, the minimum window size required for the analysis is $8\times 8$, so we choose that size here. In Figure \ref{Stokes_window}, we show the periodic extension of this window with a selected block of unknowns in blue corresponding to the upper-left pressure unknown on the computational domain.

In the following numerical results, we fix the Lam\'e parameters in such a way that the Young's modulus is $E=3\times 10^4$, and the Poisson's ratio is $\nu =0.2$. We consider different values for permeability, $K$, ranging from $K=1$ to $K=10^{-15}$. First, we consider the multigrid algorithm based on additive Schwarz smoothers using the natural weights, i.e., $w=1$ for the pressure and $w=1/9$, $1/6$, or $w=1/4$ for the displacement depending on if the degree of freedom corresponds to a vertex, to an edge, or a cell node, respectively. In Table \ref{poro_factors2}, we show the asymptotic convergence factors provided by our analysis together with the numerical ones provided by our implementation using W-cycles. We consider different numbers of smoothing steps and different values of permeability.
\begin{table}[htb]
	\begin{center}
		\begin{tabular}{ccccccccc}
			\cline{2-9}
			& \multicolumn{2}{c}{$W(1,0)$} & \multicolumn{2}{c}{$W(1,1)$} & \multicolumn{2}{c}{$W(2,1)$} & \multicolumn{2}{c}{$W(2,2)$}\\
			\cline{2-9}
			$K$ & $\rho_{2g}$ & $\rho_h$ & $\rho_{2g}$ & $\rho_h$ & $\rho_{2g}$ & $\rho_h$ & $\rho_{2g}$ & $\rho_h$\\
			\hline
			$1$        & 0.49 & 0.49 & 0.25 & 0.22 & 0.12 & 0.11 & 0.06 & 0.06 \\
			$10^{-3}$  & 0.49 & 0.49 & 0.25 & 0.21 & 0.12 & 0.11 & 0.06 & 0.06 \\
			$10^{-6}$  & 0.49 & 0.49 & 0.25 & 0.21 & 0.12 & 0.11 & 0.06 & 0.06 \\
			$10^{-9}$  & 0.65 & 0.65 & 0.42 & 0.42 & 0.27 & 0.27 & 0.18 & 0.18 \\
			$10^{-12}$ & 0.72 & 0.72 & 0.52 & 0.52 & 0.52 & 0.56 & 0.28 & 0.28 \\
			$10^{-15}$ & 0.72 & 0.72 & 0.52 & 0.52 & 0.52 & 0.56 & 0.28 & 0.28 \\
			\hline
		\end{tabular}
		\caption{Taylor-Hood discretizations of the 2D quasi-static Biot's model. Two-grid asymptotic convergence factors provided by LFA, $\rho_{2g}$, together with the computed asymptotic convergence factors, $\rho_{h}$, using several permeability values and different numbers of smoothing steps of the $51$-point additive Schwarz relaxation with natural  weights.}
		\label{poro_factors2}
	\end{center}
\end{table}

In Table \ref{poro_factors2}, we observe that there is a good match between the factors provided by LFA and the ones experimentally obtained. We also see that the performance of the multigrid method deteriorates when $K$ tends to zero, yielding a nonrobust algorithm. However, an additional advantage of our analysis is that we are able to find the optimal weights for a fixed permeability in order to obtain a robust algorithm. Since permeability is usually heterogeneous in real applications, we choose the optimal weights corresponding to the worst case scenario ($K =0$), that is, $w=0.09$ for the displacements at vertices and edges, $w=0.22$ for the displacements at interior points, and $w=1.02$ for the pressure variables.
In Table \ref{poro_opt_general}, we show the asymptotic convergence factors applying the additive Schwarz smoother with these optimal weights provided by our analysis, considering different number of smoothing steps and several permeability values.
\begin{table}[htb]
	\begin{center}
		\begin{tabular}{ccccc}
			\hline
			$K$ & \multicolumn{1}{c}{$W(1,0)$} & \multicolumn{1}{c}{$W(1,1)$} & \multicolumn{1}{c}{$W(2,1)$} & \multicolumn{1}{c}{$W(2,2)$}\\
			\hline
			$1$        & 0.58 & 0.34 & 0.19 & 0.11 \\
			$10^{-3}$  & 0.58 & 0.34 & 0.19 & 0.11 \\
			$10^{-6}$  & 0.58 & 0.34 & 0.19 & 0.11 \\
			$10^{-9}$  & 0.58 & 0.34 & 0.19 & 0.11 \\
			$10^{-12}$ & 0.60 & 0.36 & 0.21 & 0.13 \\
			$10^{-15}$ & 0.60 & 0.36 & 0.21 & 0.13 \\
			\hline
		\end{tabular}
		\caption{Taylor-Hood discretizations of the 2D quasi-static Biot's model. Two-grid asymptotic convergence factors provided by LFA, using several permeability values and different numbers of smoothing steps of the $51$-point additive Schwarz relaxation with with weights $(0.09,0.22,1.02)$.}
		\label{poro_opt_general}
	\end{center}
\end{table}

Given that $V$-cycles are less expensive than $W$-cycles, we would like to compare their performance in number of iterations. Hence, in Table \ref{poro_VW}, we show the number of iterations, $n_{iter}$, used in our implementation required to reduce the initial residual by a factor of $tol=10^{-10}$, using the optimal weights with $V$- and $W$-cycles for several smoothing steps and permeabilities.
\begin{table}[htb]
	\begin{center}
		\begin{tabular}{ccccccccc} 
			\cline{2-9}
			&  \multicolumn{8}{c}{Smoothing steps $(\nu_1, \nu_2)$}\\ 
			\cline{2-9}
			& \multicolumn{2}{c}{$(1,0)$} & \multicolumn{2}{c}{$(1,1)$} & \multicolumn{2}{c}{$(2,1)$} & \multicolumn{2}{c}{$(2,2)$}\\
			\hline
			$K$ & $V$ & $W$ & $V$ & $W$  & $V$ & $W$  & $V$ & $W$ \\
			\hline
			$1$        & 39 & 40 & 20 & 20 & 13 & 14 & 10 & 10 \\
			$10^{-3}$  & 40 & 40 & 20 & 20 & 14 & 13 & 10 & 10 \\
			$10^{-6}$  & 39 & 40 & 20 & 20 & 13 & 14 & 10 & 10 \\
			$10^{-9}$  & 40 & 40 & 23 & 20 & 14 & 13 & 10 & 10 \\
			$10^{-12}$ & 50 & 43 & 40 & 22 & 18 & 15 & 12 & 12 \\
			$10^{-15}$ & 50 & 45 & 40 & 23 & 18 & 15 & 12 & 12 \\
			\hline
		\end{tabular}
		\caption{Taylor-Hood discretizations of the 2D quasi-static Biot's model.  Number of multigrid iterations using the additive Schwarz smoothers required to reduce the initial residual by a factor of $tol=10^{-10}$, using the weights $(0.09,0.22,1.02)$ with $V$- and $W$-cycles for several smoothing steps and permeabilities.}
		\label{poro_VW}
	\end{center}
\end{table}

We conclude that $V$-cycles are also robust with respect to $K$ when applying $V(2,1)$- or $V(2,2)$-cycles since both cycles yield a similar number of iterations in order to reach the stopping criterion. Thus, the use of $V(2,2)$-cycles with the $51$-point additive Schwarz smoother applying the weights $(0.09,0.22,1.02)$ is a good approach for solving the two-field formulation of Biot's equations. 

\section{Conclusions}
\setcounter{section}{5}
\label{sec:conclusions}
A general framework for the analysis of multigrid methods based on additive Schwarz relaxations is presented in this work. This approach is based on a local Fourier analysis, which considers a basis of the Fourier space different from the classical one given in terms of the Fourier modes. The proposed analysis allows the study of a class of additive Schwarz smoothers for any problem and discretization. As examples, high-order finite-element discretizations of Poisson's problem and saddle point problems are considered since these are two formulations that benefit from applying these relaxation procedures. In addition, both additive and restricted additive Schwarz smoothers are considered. In all cases, we show an excellent match between the convergence rates provided by the analysis and those asymptotic convergence factors computationally obtained from a multigrid implementation. In the case of high-order finite-element discretizations of the Poisson problem, for each polynomial degree $p$, we find that the restricted additive Schwarz smoother achieves comparable convergence rates to the additive Schwarz smoother, having better favorable properties in terms of scalability and applicability to high-performance computing. Thus, we conclude that the restricted additive Schwarz method can be a good alternative as a relaxation procedure in a multigrid method for solving high-order discretizations of elliptic equations. Regarding the application to saddle-point problems, we consider Biot's model for poroelasticity as a model problem and show that the proposed analysis allows one to obtain optimal weights to define a multigrid solver based on additive Schwarz smoothers, which is robust with respect to physical parameters involved in the model. The proposed analysis, however, is a general approach that can be applied to any type of discretization and problem.

\section*{Acknowledgements}
\'Alvaro P\'e and Francisco J. Gaspar have received funding from the Spanish project PID2019-105574GB-I00 (MCIU/AEI/FEDER, UE). The work of Carmen Rodrigo is supported in part by the Spanish project PGC2018-099536-A-I00 (MCIU/AEI/FEDER, UE). \'Alvaro P\'e, Carmen Rodrigo and Francisco J. Gaspar acknowledge support from the Diputaci\'on General de Arag\'on, Spain (Grupo de referencia APEDIF, ref. E24\_17R). 
The work of James Adler and Xiaozhe Hu is partially supported by the National Science Foundation (NSF) under grant DMS-2208267. The research of Ludmil Zikatanov is supported in
part by the U. S.-Norway Fulbright Foundation and the U. S. National Science Foundation grant DMS-2208249.


\bibliography{mybibfile}

\end{document}